\newcommand{\N}{\mathbb{N}}               

\newcommand{\R}{\mathbb{R}}
\newcommand{\Z}{\mathbb{Z}}
\newcommand{\A}{\mathscr{A}}
\newcommand{\Sp}{\mathbb{S}}
  
\newcommand{\m}{\mathfrak{m}}  
\newcommand{\Diff}{{\rm Diff}}   
\newcommand{\E}{{\mathscr E}}   
\mathchardef\varepsilon="010F
\mathchardef\epsilon="0122
\mathchardef\vartheta="0112
\mathchardef\theta="0123
\mathchardef\varrho="011A
\mathchardef\rho="0125
\mathchardef\varphi="011E     
\mathchardef\phi="0127
\renewcommand \emptyset \varnothing
\documentclass[12pt]{article}            
\usepackage[all]{xy}
\usepackage[latin1]{inputenc}        
\usepackage[dvips]{graphics}
\usepackage[dvips]{graphicx}
\usepackage{amsfonts}
\usepackage{amssymb}
\usepackage{amsmath}
\usepackage{amstext}
\usepackage{amsbsy}
\usepackage{amsopn}
\usepackage{amsthm}
\usepackage{amscd}
\usepackage{amsxtra}
\usepackage{mathrsfs}
\usepackage{upref}
\author{\sc Gianmarco Capitanio\footnote{Research partially supported by Istituto Nazionale Di Alta Matematica
F. Severi.}}
\date{}
\title{Stable tangential families and singularities of their envelopes}

\begin{document}        

\numberwithin{equation}{section}                
\theoremstyle{plain}
\newtheorem{theorem*}{\bf Theorem}
\newtheorem{theorem}{\bf Theorem}
\newtheorem{lemma}{\bf Lemma}
\newtheorem{step}{\bf Step}
\newtheorem{proposition}{\bf Proposition}
\newtheorem*{proposition*}{\bf Proposition}
\newtheorem{corollary}{\bf Corollary}
\newtheorem*{corollary*}{\bf Corollary}
\theoremstyle{definition}
\newtheorem*{definition*}{\bf Definition}
\newtheorem*{definitions*}{\bf Definitions}
\newtheorem*{conjecture}{\bf Conjecture}
\newtheorem{example}{\bf Example}
\newtheorem*{example*}{\bf Example}
\theoremstyle{remark}
\newtheorem*{remark*}{\bf Remark}
\newtheorem{remark}{\bf Remark}
\newtheorem*{acknowledgements}{\bf Acknowledgements}

\maketitle
\begin{abstract}
We study tangential families, i.e. systems of rays emanating
tangentially from given  curves. 
We classify, up to Left-Right equivalence, stable singularities of
tangential family germs (under deformations among tangential families)
and we study their envelopes.  
We discuss applications of our results to the case of tangent
geodesics of a curve. 
\end{abstract}


\section{Introduction}\label{sct:1}

A tangential family is a system of rays emanating tangentially from a
support curve. 
Tangential families naturally arise in Geometry of Caustics (see
\cite{arnold2001}) and in Differential Geometry. 
For example, the tangent geodesics of a curve in a Riemannian
surface define a tangential family.
Envelopes of tangential families with singular support are
studied in  \cite{io2002} and \cite{io200?}. 

The roots of the theme go back to Huygens' investigation  of caustics
of rays of light. 
The theory of tangential families is a part of Envelope theory. 
Thom showed in \cite{Thom} that the singularities of envelopes of
generic $1$-parameter families of smooth plane curves are
semicubic cusps and transversal self-intersections.
Normal forms of generic families of plane curves near regular
points of the envelope have been found by V.I. Arnold (see
\cite{arnold1976a} and \cite{arnold1976b}).  

However, our study differs from the approaches of Thom and
Arnold. 
For example, for the first time the situation when the support curve
is not the only local component of the envelope is studied. 

The aim of this paper is to classify, up to Left-Right equivalence, 
the singularities of tangential family germs which are stable under small
deformations among tangential families, and to study their envelopes.  
We prove that in addition to a regular envelope, there exists just one
more local stable singularity of envelopes, the second order
self-tangency. 
Applications to families of geodesics emanating tangentially from a
regular curve on a Riemannian surface are given. 

\begin{acknowledgements}
I wish to express my deep gratitude to V.I. Arnold for suggesting me to
study this problem and for his careful reading of this paper.
I also like to thank M. Garay for very useful discussions and
comments, M. Kra\"iem and F. Aicardi for their encouragement.
\end{acknowledgements}

\section{Preliminary definitions}\label{sct:2}

Unless otherwise specified, all the objects considered below are 
supposed to be of class $\mathscr C^\infty$; by {\it curve} we mean an 
embedded $1$-submanifold of $\mathbb{R}^2$. 

Let us consider a map $f:\R^2\rightarrow \R^2$ of the plane $\R^2$,
whose coordinates are denoted by $\xi$ and $t$.  
If $\partial_t f$ vanishes nowhere, then $f$ defines the $1$-parameter
family of the plane curves parameterized by $f_\xi:=f(\xi,\cdot)$; 
these curves may have double points.   
The map $f$ is called a {\it parameterization} of the family.

\begin{definition*}
The family parameterized by $f$ is called a {\it tangential family} if 
$\partial_\xi f$ and $\partial_t f$ are parallel non zero vectors at
every point $(\xi,t=0)$, and the image of $f(\cdot,0)$ is an 
embedded curve, called the {\it support} of the family. 
\end{definition*}

We remark that the curves parameterized by $f_\xi$ are tangent to 
the support $\gamma$ at $f(\xi,0)$. 
The {\it graph} of the tangential family is the surface 
\[\Phi:= \left\{ (q,p) : q=f(\xi,0),\  p=f(\xi,t) ,\  
\xi,t\in\R \right\} \subset \gamma\times \mathbb{R}^2 \ . \]
Let us consider the two natural projections of $\Phi$ on $\gamma$ and
$\R^2$, $\pi_1 : (q,p)\mapsto q$ and $\pi_2 : (q,p)\mapsto p$. 
The first projection $\pi_1$ is a fibration; the images by $\pi_2$ of
its fibers are the curves of the family.

The {\it criminant set} of the tangential family is the critical set
of ${\pi_2}$; the {\it envelope} is the apparent 
contour of its graph in the plane (i.e., the critical 
value set of $\pi_2$). 
By the very definition, the support of a tangential family belongs to
its envelope.  

Our study of tangential families being local, we consider their 
parameterizations as elements of $(\m_{\xi,t})^2$, where $\m_{\xi,t}$ is the
space of function germs in two variables vanishing at the origin. 
We denote by $X_0$ the subset of $(\m_{\xi,t})^2$ formed by all the map
germs parameterizing a tangential family.

\begin{remark*} 
The graph of a tangential family germ is smooth.
\end{remark*}

Let us denote by $\Diff (\R^2,0)$ the group of the
diffeomorphism germs of the plane keeping fixed the origin, and by
$\A$ the direct product $\Diff (\R^2,0)\times \Diff(\R^2,0)$.
Then the group $\A$ acts on $(\m_{\xi,t})^2$ by the rule 
$(\phi,\psi) \cdot f := \psi\circ f\circ\phi^{-1}$.
Two map germs of $(\mathfrak m_{\xi,t})^2$ are said to be 
{\it Left-Right equivalent}, or {\it $\A$-equivalent}, 
if they belong to the same $\mathscr A$-orbit. 
The {\it singularity} of a tangential family germ is its
$\A$-equivalence class.  
Note that $X_0$ is not $\A$-invariant.

\begin{remark*}
Critical value sets of $\A$-equivalent map germs are
diffeomorphic; in particular, envelopes of $\A$-equivalent
tangential families are diffeomorphic.
\end{remark*}

In some situations, as for instance in the study of geodesic tangential
family evolution under small perturbations of the metric, it would be natural 
to perturb a tangential family only among tangential families.

\begin{definition*} 
A $p$-parameter {\it tangential deformation} of a tangential family
$f:\R^2\rightarrow \R^2$ is a mapping 
$F:\R^2\times \R^p \rightarrow \R^2$, 
such that $F_\lambda:=F(\cdot;\lambda)$ is a
tangential family for every $\lambda$ and $F_0=f$.
\end{definition*} 

For example, the translation of the origin is a $2$-parameter
tangential deformation.  
An analoguous definition holds for germs.
The notion of tangential deformation of a tangential family germ can
be extended to its whole $\A$-orbit via the action of the group $\A$
on $\m_{\xi,t}$. 

\begin{remark*}
A tangential deformation induces a smooth deformation on the
family support.  
\end{remark*}

A singularity is said to be {\it stable} if,  
for every representative $f$ of it, every tangential deformation 
$F_\lambda$ of $f$ has a singularity $\A$-equivalent to that of $f$ 
at some $\lambda$-depending point, arbitrary close to the origin 
for $\lambda$ small enough. 

We end this section recalling some definitions 
(see e.g. \cite{ST1}, \cite{Martinet}). 
The {\it extended tangent space} at $f\in(\m_{\xi,t})^2$ to its
$\A$-orbit is the subspace of $\mathscr E^2_{\xi,t}$ defined by 
\[ T_e\A(f):= \langle \partial_\xi f , \partial_t f\rangle_{\mathscr
  E_{\xi,t} }  +  f^*(\mathscr E_{x,y})\cdot  \mathbb{R}^2  \ , \] 
where $\E_{\xi,t}$ is the ring of the function germs
$(\R^2,0)\rightarrow\R$ in the variables $\xi,t$ and the homomorphism 
$f^*:\mathscr E_{x,y}\rightarrow \mathscr
E_{\xi,t}$ is defined by $f^*g:=g\circ f$.  
Every $\mathscr E_{\xi,t}$-module can be viewed as an
$\mathscr E_{x,y}$-module via this homomorphism.
The {\it extended codimension} of $f$ is the dimension of the
quotient space $\mathscr E^2_{\xi,t}/T_e\A(f)$ as a real vector space. 

\section{Stable tangential family germs}\label{sct:3}

In this section we state our main results.  
The proofs are given in section \ref{sct:5}.

\begin{theorem}\label{thm:1}
All the stable singularities of tangential family germs are those
listed in the table below, together with their extended
codimension. 
\begin{center}
\begin{tabular}{|| c | c | c ||} 
\hline \hline
Singularity & Representative & codim\\ \hline\hline
$\rm I$& $(\xi+t,t^2)$& $0$ \\ \hline 
${\rm II}$ & $(\xi+t,t^2\xi)$&$1$ \\ \hline \hline
\end{tabular}
\end{center}
\bigskip
\end{theorem}

\begin{corollary*}
The envelope of any tangential family germ having a singularity 
${\rm I}$ is smooth, while the envelope of any tangential family germ
having a singularity ${\rm II}$ has an order $2$ self-tangency.  
\end{corollary*}

\begin{figure}[htpb]
  \centering
   \scalebox{.35}{\input{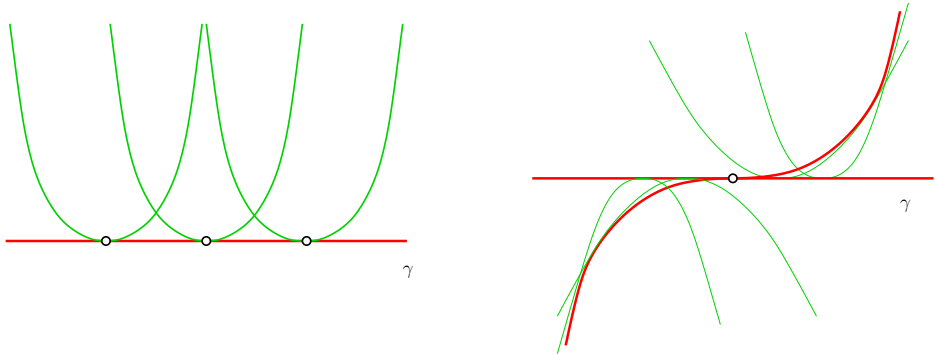}}
  \caption{Stable tangential family germs and their envelopes.}
 \label{stablefig1}
\end{figure}

Tangential families representing stable singularities 
${\rm I}$ and ${\rm II}$, together with their envelopes, are
depicted in figure \ref{stablefig1}.

\begin{remark*}
Singularity ${\rm I}$ is a fold. 
Its extended codimension being $0$, every small enough deformation of
it is tangential. 
On the other hand, the extended codimension of singularity 
${\rm II}$ is equal to $1$: there exist non tangential deformations.  
Under such deformations, the family envelope experiences a
beak to beak perestroika. 
\end{remark*}

Let us consider a tangential family germ $f\in X_0$.
The fiber $\pi_2^{-1}(0,0)$ of the projection $\pi_2:\Phi\rightarrow
\R^2$ defines a {\it vertical direction} in the tangent plane of the
family graph $\Phi$ at the origin.
Note that the branch of the criminant set of $f$, 
projecting into the support, is non vertical.

\begin{definition*}
A tangential family germ is said to be of {\it first type} if 
its criminant set has only one branch, of {\it second type} if it has exactly 
two branches, provided that they are smooth, transversal and
non vertical. 
\end{definition*}
 
Stable tangential family germs admit the following
characterization. 

\begin{theorem}\label{thm:2}
A tangential family germ has a singularity ${\rm I}$ 
(resp., ${\rm II}$) if and only if it is of first type (resp., of
second type). 
\end{theorem}
 
\section{Examples and applications}\label{sct:4}

The simplest example of tangential family is that of the tangent lines
to a curve in the Euclidean plane.   
This example can be generalized replacing the Euclidean structure 
by a Riemannian structure: every curve in a Riemannian
surface $(M,g)$ defines the {\it geodesic tangential family} formed by its
tangent geodesics.  
By simplicity, we assume $M$ compact. 

Generically, a smooth map $f:\Sp^1\rightarrow M$ defines a curve, 
whose singularities may be only transversal self-intersections.     
In this case, we call {\it geodesic envelope of $f$} the envelope of
the geodesic tangential family of support $f(\Sp^1)$.
Here, as well as in Theorem \ref{thm:3} below, ``generically'' means 
that the mappings satisfying the claim form an open dense subset of
$C^\infty(\Sp^1,M)$ for the Whitney topology.

A deformation of the metric induces a tangential deformation on the 
geodesic envelope of $f$ (fixing the support); similarly, a small
enough deformation of the function $f$ induces a tangential deformation of
its geodesic envelope. 

\begin{theorem}\label{thm:3}
The geodesic envelope of a smooth map $f:\Sp^1\rightarrow M$ 
is generically a curve, whose singularities may be only 
transversal self-intersections, semicubic cusps and order $2$
self-tangencies. 
All the envelope singular points form a discrete set; 
the envelope order $2$ self-tangencies coincide with 
the inflection points of $f(\Sp^1)$. 
Moreover, this envelope is stable under small perturbations of $f$ and
of the metric.
\end{theorem}

Thorem \ref{thm:3} follows from Theorem \ref{thm:1} and Thom's
Transversality Lemma.
A similar statement holds if we replace the Riemannian structure with a 
projective structure. 

We discuss now an example of a global geodesic tangential
family whose envelope has infinitely many branches.

Let us consider the sphere $\Sp^2\subset\R^3$ equipped with
the induced metric.  
Let $\gamma_r\subset \Sp^2$ be a circle of radius $r$ in $\R^3$ and
let $f:\Sp^1\times \R\rightarrow \Sp^2$ be a parameterization of the
geodesic tangential family of support $\gamma_r$ in
$\Sp^2$, where for every $\xi\in \Sp^1\simeq \gamma_r$, the map 
$f_\xi:\R\rightarrow \Sp^2$ is a $2\pi$-periodical
parameterization of the $\Sp^2$ great circle at $\xi$. 

The family envelope has two branches, $\gamma_r$ and its opposite
circle $-\gamma_r$, each one having infinite multiplicity. 
Indeed, the curve $\det Df(\xi,t)=0$ has infinitely many branches 
$$C_n:=(\xi\in S^1, t=n\pi) \ , \quad n\in \Z\ .$$ 
Note that $E_n=(-1)^n\gamma_r$, where $E_n:= f(C_n)$ is called the
{\it order $n$ envelope} of $\gamma_r$. 

Let $\tilde \Sp^2$ be a small perturbation of the sphere.  
Denote by $\tilde \gamma_r$ the support of the perturbed
family and by $\tilde f:\Sp^1\times\R\rightarrow \tilde \Sp^2$ 
a parameterization of the perturbed geodesic tangential family. 
Fix $N,M\in \N$ arbitrary large, $N<M$.  
 
\begin{theorem}\label{thm:4}
If $\tilde \Sp^2$ is close enough to $\Sp^2$, then the order $n$ 
envelopes of $\tilde \gamma_r$ are generic spherical 
closed caustics with zero Maslov number for
$|n|<M$, smooth for $|n|<N$. 
\end{theorem}

\begin{proof}
If $\tilde \Sp^2$ is close enough to $\Sp^2$, for every $n<M$
the curves $\tilde C_n$ are pairwise disjoint small perturbations of
the curves $C_n$.  
Each order $n$ envelope is a small perturbation of the
corresponding order $n$ unperturbed envelope. 
Since the unperturbed tangential family is of
first type at every point of its envelope, by Theorem \ref{thm:1} 
the perturbed family is of first type at any point of the envelope
branches $\tilde E_n$ for every $|n|<N$, provided that $\tilde \Sp^2$
is close enough to $\Sp^2$. 
The perturbed envelope $\tilde E_n$ can be viewed
as a caustic, images of $\tilde E_{n-1}$ under a geodesic flow.  
Hence, generic singularities of caustics, as semicubic cusps and transversal 
self-intersection, may arise in envelopes $\tilde E_n$ for  $N<n<M$. 
\end{proof}

Consider now a fixed arbitrary small perturbation $\tilde \Sp^2$ of
$\Sp^2$.  

\begin{theorem}\label{thm:5}
The first order envelope of $\tilde
\gamma_r$ has at least four cusps, provided that $\tilde \Sp^2$ is a
generic convex surface close to the sphere and $r$ is small enough.  
\end{theorem}

\begin{proof}
When $r\rightarrow 0$, $\tilde \gamma_r$ shrinks about a point of
$\tilde S^2$ and the envelope branches $\tilde E_n$, defined for
$n<M$, approach the order $n$  caustics of this point. 
If $\tilde \Sp^2$ is a generic convex surface close
enough to $\Sp^2$, then the Last Geometrical Theorem of Jacobi states that the
first caustic of any generic point has at least $4$ cusps (see
\cite{arnold1999a}).   
\end{proof}

For a fixed $n$, the caustic of order $n$ of a point
has at least four cusps, provided that the perturbation is small
enough. 
But whether it holds simultaneously for all the values $n$ provided
that the perturbation is small enough is still a conjecture (see
\cite{arnold1999a}).  
This conjecture can be generalized to envelopes of closed convex smooth
curves of small length on convex surfaces. 

\section{Proof of Theorem \ref{thm:1} and \ref{thm:2}}\label{sct:5}

In order to prove Theorem \ref{thm:1} and \ref{thm:2}, we introduce now a 
{\it prenormal form} to which is possible to bring every tangential 
family germ by an $\A$-equivalence.  
This prenormal form does not belong to $X_0$ and it is not unique. 
We denote by $\delta(n)$ any function in two variables with 
vanishing $n$-jet at the origin.

\begin{lemma}\label{lemma:1}
Every germ $f\in X_0$ is $\mathscr A$-equivalent to
a map germ of the form 
\[ (\xi,t)\mapsto\left(\xi, k_0 t^2 + (\alpha-k_1) t^3+k_1 t^2\xi+t^2 \cdot 
\delta(1)\right)\ , \] 
for some coefficients $k_0,k_1,\alpha\in\R$.
\end{lemma}

\begin{proof} 
Let us consider a representative of the tangential family germ,
defined in an arbitrary small neighborhood $\cal U$ of the origin. 
Fix a new coordinate system $\{x,y\}$ in $\mathcal U$, such that 
the support is $y=0$.  
Denote by $\Gamma_\xi$ the family curve
corresponding to the support point $(\xi,0)$.    
Define $k_0$ and $k_1$ by the expansion $k_0 + k_1 \xi + o(\xi)$ 
(for $\xi\rightarrow 0$) of half of the curvature of $\Gamma_\xi$ at
$(\xi,0)$;   
$\alpha$ is similarly defined by the expansion 
$k_0 t^2+\alpha t^3+o(t^3)$ for $t\rightarrow 0$ 
of the function, whose graph in $\mathcal U$ is $\Gamma_0$. 
For any small enough value of $\xi$, $\Gamma_\xi$ can be 
parameterized near $(\xi,0)$ by its projection $\xi+t$ on the $x$-axis. 
In this manner we get a parameterization of the family, that can be
written as 
$$(\xi,t)\mapsto \big(t+\xi,k_0t^2+k_1t^2\xi+\alpha t^3+t^2\cdot 
\delta(1)\big)\ .$$ 
Taking $\xi+t$ and $t$ as new parameters, we obtain the
required germ. 
\end{proof}

The type of a tangential family germ is determined by its
prenormal form. 

\begin{lemma}\label{lemma:2}
A tangential family germ (in prenormal form) 
is of first type if and only if $k_0\not=0$;
it is of second type if and only if $k_0=0$ and 
$k_1\not= 0,\alpha$.
\end{lemma}

Lemma \ref{lemma:2} follows from Lemma \ref{lemma:1} by a direct computation. 

\begin{remark*}
Consider any representative $\bar f$ of $f\in X_0$.  
Then $f$ is of first type if and only if the curve
$\Gamma_\xi$, parameterized by $\bar f_\xi$, has an order $1$ tangency 
with the support $\gamma$ at $\bar f(\xi,0)$ for every small enough $\xi$.  
Moreover, if $f$ is of second type, then the tangency
order of $\gamma$ and $\Gamma_0$ is at least $2$, the tangency order
of $\Gamma_\xi$ and $\gamma$ being $1$ for every small enough $\xi\not=0$.
\end{remark*}

The proof of Theorem \ref{thm:2} is divided into two steps. 

\begin{step}\label{step:1}
Every tangential family germ of first type is $\A$-equivalent to the
representative $f_{\rm I}(\xi,t):=(\xi+t,t^2)$ of singularity ${\rm I}$.  
\end{step}

\begin{proof} 
By Lemmas \ref{lemma:1} and \ref{lemma:2}, every first type
tangential family germ is $\mathscr A$-equivalent to a map germ of the
form $(\xi,t^2+\delta (2))$.  
Now, it is well known that the $2$-jet $(\xi,t^2)$ is $\A$-sufficient
and stable, so its extended codimension is $0$.
Finally, $(\xi,t^2)$ is $\A$-equivalent to $f_{\rm I}$.
\end{proof}

\begin{step}\label{step:2}
Every tangential family germ of second type is $\A$-equivalent to the
representative $f_{\rm II}(\xi,t):=(\xi+t,\xi t^2)$ of singularity ${\rm II}$.  
\end{step}

\begin{proof} 
Let us set $h(\xi,t):=(\xi,\xi t^2+t^3)$; this germ is
$\A$-equivalent to $f_{\rm II}$. 
By Lemmas \ref{lemma:1} and \ref{lemma:2}, 
every second type tangential family germ is $\A$-equivalent to a germ 
of the form $h+\big(0,P_n+\delta(n)\big)$, 
for some homogeneous polynomial $P_n$ of degree $n>3$. 
We prove now that every such a germ is 
$\mathscr A$-equivalent to a germ $h+\big(0,\delta(n)\big)$. 
This equivalence provides by induction the formal 
$\mathscr A$-equivalence between the initial prenormal form and $h$. 

For any fixed $i\in\{2,\dots,n\}$, we can kill the term on
$\xi^{n-i}t^{i}$ in $P_n$, changing the coefficient of 
$\xi^{n-i+1} t^{i-1}$ and higher order terms. 
This is done by the coordinate change 
$$(\xi,t)\mapsto
(\xi,t+a\xi^{n-i}t^{i-2})\ , $$ for a suitable $a\in\R$. 
Hence, we may assume
$P_n(\xi,t)=A\xi^{n-1}t+B\xi^n$.  
Moreover, we can also suppose $B=0$, up to the coordinate change
$(x,y)\mapsto (x,y-Bx^n)$. 
Now, the coordinate change 
$(\xi,t)\mapsto(\xi+3A\xi^{n-2}/2,t-A\xi^{n-2}/2)$ takes our germ to
the form $h+(A\xi^{n-2},\delta(n))$, which is $\mathscr
A$-equivalent to $h+\big(0,\delta(n)\big)$.

We compute now the extended codimension of $h$. 
Since $\langle h\rangle _{\mathscr{E}_{\xi,t}} =
\langle \xi,t^3\rangle_{\mathscr{E}_{\xi,t}}$, 
the $\mathscr E_{x,y}$-module $\E_{\xi,t}$ is
generated by $1$, $t$ and $t^2$ (Preparation Theorem).  
Now it is easy to check that all the vector monomials $(\xi^pt^q,0)$
for $p+q\geq 0$ and $(0,\xi^pt^q)$ for $p+q\geq 2$ belong to the
tangent space $T_e\A(h)$, so we have
\begin{equation}\label{eq:1}
\mathscr E^2_{\xi,t} = 
T_e\A(h) + \E_{x,y}\cdot \begin{pmatrix} 0\\t\end{pmatrix} = 
T_e\A(h) \oplus \R\cdot  \begin{pmatrix} 0\\t\end{pmatrix} ;
\end{equation} 
indeed, $(0,t)$ does not belong to $T_e\A(h)$. 
This proves that the extended codimension of $h$ is equal to $1$. 
In particular, $h$ is finitely determined.
\end{proof}

We prove now Theorem \ref{thm:1}. 
The fold $f_{\rm I}$ is stable under any small deformation (see
e.g. \cite{Martinet}), in particular under small
tangential deformations.    
Equality (\ref{eq:1}) implies that the mapping 
$H:\R^2\times\R \rightarrow\R^2$, defined by 
$H(\xi,t;\lambda):= h(\xi,t)+\lambda (0,t)$, is an $\mathscr
A$-miniversal deformation of $h$. 
Under this deformation, the critical set of $h$ experiences a beak
to beak perestroika, so the deformation is not tangential. 

This allows us to prove that singularity ${\rm II}$ is stable under small
tangential deformations. 
Let us consider a $p$-parameter tangential deformation $K$ of $h$. 
By the very definition of versality, $K$ can be represented as 
$$K(\xi,t;\mu)\equiv \Psi(H(\Phi(\xi,t;\mu),\Lambda(\mu)),\mu)\ ,$$
where $\Phi(\cdot;\mu)$ and $\Psi(\cdot;\mu)$ are deformations
of the identity diffeomorphisms of the source and the target planes 
$\R^2$ and $\Lambda$ is a function germ $(\R^p,0)\rightarrow(\R,0)$. 
Therefore, for any fixed $\mu$, the critical value sets of
$K(\cdot,\mu)$ and $H(\Phi(\cdot;\mu),\Lambda(\mu))$ are
diffeomorphic. 
Assume $\Lambda\not\equiv 0$; then there exists an arbitrary
small $\mu_0$ for which  $\Lambda(\mu_0)\not=0$. 
Hence, when $\mu$ goes from $0$ to $\mu_0$, the envelope of the
deformed family $K_\mu$ experience a beak to beak perestroika.  
Since the deformation $K$ is tangential, this is impossible, so
$\Lambda\equiv 0$ and $H$ is trivial. 
This proves that singularity ${\rm II}$ is stable.  

Let us consider now a tangential family germ $f$ in prenormal form.   
Then the $1$-parameter tangential deformation
$f(\xi,t)+\lambda(0,t^3+\xi t^2)$ is non trivial whenever $f$ is
neither of first nor second type.   
Indeed, for any small enough $\lambda$, the deformed family has a
singularity of type ${\rm II}$ at the origin.  
This proves that every tangential family germ, stable under small
tangential deformations, has a singularity ${\rm I}$ or ${\rm II}$.
The proof of Theorem \ref{thm:1} is now completed.


\bigskip

\noindent {\it Gianmarco Capitanio} \\
\noindent {\it    Universit{\' e} D. Diderot -- Paris VII} \\ 
\noindent {\it    UFR de Math{\' e}matiques} \\
\noindent {\it    Equipe de G{\' e}om{\' e}trie et Dynamique}  \\
\noindent {\it    Case 7012 -- 2, place Jussieu} \\
\noindent {\it    75251 Paris Cedex 05} \\
{\it  e-mail:} Gianmarco.Capitanio@math.jussieu.fr 
  

\addcontentsline{toc}{section}{Bibliographie}

\begin{thebibliography}{abcd}

\bibitem{arnold1976a} {\sc Arnold, V. I.}: On the envelope theory,
  {\it Uspecki Math. Nauk.} 3, {\bf 31} (1976), 248--249 (in
  Russian). 

\bibitem{arnold1976b} {\sc Arnold, V. I.}: Wave front evolution and
  equivariant Morse lemma, {\it Comm. Pure Appl. Math.} 6, {\bf 29}
  (1976), 557--582. 

\bibitem{arnold1999a}
{\sc Arnold, V. I.}:  
Topological problems in wave propagation theory and topological
economy principle in algebraic geometry.  
The Arnoldfest (Toronto, ON, 1997),  39--54,  
Fields Inst. Commun., 24, Amer. Math. Soc., Providence, RI, 1999. 

\bibitem{arnold2001} {\sc Arnold, V. I.}: Astroidal geometry of
  hypocycloids and the Hessian topology of hyperbolic polynomials, 
{\it Russ. Math. Surv.} 6, {\bf 56} (2001)
  1019--1083. 

\bibitem{ST1} {\sc Arnold, V. I.; Goryunov, V. V.; Lyashko, O. V.;
Vasiliev, V. A.}: {\it Singularity theory. I.}  
Springer-Verlag, Berlin, 1998. 

\bibitem{io2002} {\sc Capitanio, G}: On the envelope of 1-parameter
  families of curves tangent to a semicubic cusp, {\it
    C. R. Math. Acad. Sci. Paris} 3, {\bf 335} (2002), 249--254. 

\bibitem{io200?} {\sc Capitanio, G.}: Singularities of envelopes of
  curves tangent to a semicubic cusp, to appear in Proc. of Suzdal
  Int. Conf. 2002.  

\bibitem{Martinet} {\sc Martinet, J.}: {\it Singularities of smooth
    functions and maps}, L. M. S. Lecture Note Series, 58. Cambridge
    University Press, 1982.  

\bibitem{Thom} {\sc Thom, R.}: Sur la théorie des enveloppes, {\it
    J. Math. Pures Appl.} 9, {\bf 41} (1962), 177--192. 

\end{thebibliography}
\end{document}